\documentclass[10pt,twoside]{article}
\usepackage{latexsym}
\usepackage{Latex-document}

\title{\bf Heat Kernels and the Index Theorems on\vskip -2mm Even and Odd Dimensional
Manifolds\thanks{Partially supported by  the MOEC  and the 973 Project.}\vskip 6mm}

\author{Weiping Zhang\vspace*{-0.5cm}\thanks{Nankai \ Institute \ of \ Mathematics, \ Nankai \ University, \ Tianjin \ 300071, \ China. \ E-mail: weiping@nankai.edu.cn}}

\date{\vspace{-8mm}}

\begin{document}

\maketitle

\thispagestyle{first} \setcounter{page}{361}

\begin{abstract}

\vskip 3mm

In this talk, we review the heat kernel approach to the Atiyah-Singer index theorem for Dirac operators on closed
manifolds, as well as the Atiyah-Patodi-Singer index theorem for Dirac operators on manifolds with boundary. We
also discuss the odd dimensional counterparts of the above results. In particular, we describe a joint result with
Xianzhe Dai on an index theorem for Toeplitz operators on odd dimensional manifolds with boundary.

\vskip 4.5mm

\noindent {\bf 2000 Mathematics Subject Classification:} 58G.

\noindent {\bf Keywords and Phrases:} Index theorems, heat kernels, eta-invariants, Toeplitz operators.
\end{abstract}

\vskip 12mm

\section{Introduction}

\vskip-5mm \hspace{5mm}

As is well-known, the  index theorem proved by Atiyah and Singer [AS1]
in 1963, which expresses the analytically defined index of elliptic differential
operators through purely topological terms, has had a wide range of
implications in mathematics as well as in mathematical physics. Moreover, there have
been up to now many different proofs of this celebrated result.

The existing proofs of the Atiyah-Singer index theorem can roughly be divided into three
categories:

(i) The cobordism proof: this is the proof originally given in [AS1]. It uses the cobordism
theory developed by Thom and modifies Hirzebruch's proof of his Signature theorem as well
as his Riemann-Roch theorem;

(ii) The $K$-theoretic proof: this is the proof given  by Atiyah and Singer in [AS2]. It
modifies Grothendieck's proof of the Hirzebruch-Riemann-Roch theorem and relies on
the topological $K$-theory developed by Atiyah and Hirzebruch. The Bott periodicity theorem plays
an important role in this proof;

(iii) The heat kernel proof: this proof originates from a simple and beautiful formula
due to Mckean and Singer [MS], and has closer relations with differential geometry as well
as mathematical physics. It also lead directly to the important
Atiyah-Patodi-Singer index theorem
for Dirac operators on manifolds with boundary.

In this article, we will survey some of the developments concerning the heat kernel proofs
of various index theorems, including a recent result with Dai [DZ2]
on an index theorem for Toeplitz
operators on odd dimensional manifolds with boundary.

\section{Heat kernels and the index theorems on even dimensional manifolds}

\vskip-5mm \hspace{5mm}

We start with a smooth closed oriented $2n$-dimensional
manifold $M$ and two smooth complex
vector bundles $E$, $F$
over $M$, on which there is an elliptic differential operator between the spaces
of smooth sections,
$D_+:\Gamma(E) \rightarrow \Gamma(F) .$

If we equip $TM$ with a Riemannian metric and $E$, $F$ with Hermitian metrics repectively,
then $\Gamma(E) $ and $\Gamma(F) $ will carry canonically induced inner products.

Let
$D_-:\Gamma(F) \rightarrow \Gamma(E) $
be the formal adjoint of $D_+$ with respect to these inner products. Then the index
of $D_+$ is given by
$$ {\rm ind}\, D_+= \dim \left(\ker D_+\right)-\dim\left(\ker D_-\right).\eqno(2.1)$$
It is a topological invariant not depending on the metrics
on $TM$, $E$ and $F$.

The famous Mckean-Singer formula [MS] says that ${\rm ind}\, D_+ $ can also be computed
by using the heat operators associated to the Laplacians $D_-D_+$ and $D_+D_-$.
That is, for any $t>0$, one has
$$ {\rm ind}\, D_+ = {\rm Tr} \left[\exp\left(-tD_-D_+\right)\right]
- {\rm Tr} \left[\exp\left(-tD_+D_-\right)\right] .\eqno(2.2)$$
By introducing the ${\bf Z}_2$-graded vector bundle $E\oplus F$ and setting
$D=\mbox{\scriptsize $\left(\begin{array}{cc} 0& D_-\\ D_+ & 0\end{array}\right)$}$,
we can rewrite the difference of the two traces in the right hand side of (2.2)
as a single ``supertrace'' as follows,
$${\rm ind}\, D_+ = {\rm Tr}_s\left[ \exp\left(-tD^2\right) \right] ,
\ \ \ \ \mbox{for any $t>0$}.\eqno(2.2)'$$

Let $P_t(x,y)$ be the smooth kernel of $\exp(-tD^2) $ with respect to the volume
form on $M$. For any $f\in \Gamma(E\oplus F)$, one has
$$\exp\left(-tD^2\right) f (x)=\int_M P_t(x,y) f(y)dy.\eqno(2.3)$$
In particular,
$$ {\rm Tr}_s \left[ \exp\left(-tD^2\right) \right] =\int_M
{\rm Tr}_s \left[P_t(x,x) \right] dx.\eqno(2.4)$$

Now, for simplicity, we assume that the elliptic operator $D$ is of order one.
Then by a standard result, which goes back to Minakshisundaram and Pleijiel [MP],
one has that when $t>0$ tends to $0$,
$$ P_t(x,x) ={1\over (4\pi t)^n}
\left(a_{-n}+a_{-n+1}t +\cdots +
a_0t^n+o_x\left(t^n\right)\right),\eqno(2.5)$$
where $a_i\in {\rm End}((E\oplus F)_x)$, $i=-n,\dots, 0$.

By (2.2)', (2.4) and (2.5), and by taking $t>0$ small enough, one deduces that
$$ \int_M{\rm Tr}_s [a_i]dx =0, \ \ \ -n\leq i< 0,$$
$${\rm ind}\, D_+ =\left({1\over 4\pi}\right)^n\int_M{\rm Tr}_s [a_0]dx .\eqno(2.6)$$

Mckean and Singer conjectured in [MS] that for certain geometric
operators, there should be some ``fantastic cancellation'' so that the following
far reaching refinement of (2.6) holds,
$$ {\rm Tr}_s [a_i] =0, \ \ \ -n\leq i< 0, $$
and moreover, ${\rm Tr}_s [a_0]$ can be calculated simply in the Chern-Weil geometric
theory of characteristic classes.

In fact, as a typical example, let $M$ be an even dimensional compact smooth oriented {\it spin} manifold carrying
a Riemannian metric $g^{TM}$. Let $R^{TM}$ be the curvature of the Levi-Civita connection associated to $g^{TM}$.
Let $S(TM)=S_+(TM)$ $\oplus S_-(TM)$ be the Hermitian vector bundle of $(TM, g^{TM})$-spinors, and \linebreak
$D_+:\Gamma(S_+(TM))\rightarrow \Gamma(S_-(TM))$ the associated {\it Dirac} operator.

One then has the formula (cf. [BGV, Chap. 4, 5]),
$$ \lim_{t\rightarrow 0} {\rm Tr}_s \left[ P_t(x,x) \right]dx=
\left\{ \widehat{A} \left( {R^{TM}\over 2\pi}\right) \right\}^{\rm max} :=
\left\{ {\det}^{1/2}\left({{\sqrt{-1}\over 4\pi} R^{TM}\over
\sinh \left({\sqrt{-1}\over 4\pi }R^{TM}\right)}\right) \right\}^{\rm max},\eqno (2.7)$$
which implies the Atiyah-Singer index theorem [AS1] for $D_+$:
$${\rm ind}\, D_+=\widehat{A} (M):=\int_M
\widehat{A} \left( {R^{TM}\over 2\pi}\right) .\eqno(2.8)$$

A result of type (2.7) is called a {\it local index theorem}. The first proof of
such a local result was given by V. K. Patodi [P] for the de Rham-Hodge
operator $d+d^*$. Other direct heat kernel
proofs of (2.7) have been given by Berline-Vergne, Bismut,
Getzler and Yu respectively. We refer to [BGV] and [Yu] for more details.

The heat kernel proof of the local index theorem leads to  a generalization of the
index theorem for Dirac operators to the case of manifolds with boundary. This was
achieved by Atiyah, Patodi and Singer in [APS], and will be reviewed in the next section.

\section{The index theorem for Dirac operators
on even dimensional manifolds with boundary}

\vskip-5mm \hspace{5mm}

Let $M$ be a smooth compact
 oriented even dimensional {\it spin} manifold with (nonempty)
  smooth boundary $\partial M$.
Then $\partial M$ is again oriented and {\it spin.}

Let $g^{TM}$ be a metric on $TM$. Let $g^{T\partial M}$ be its restriction on
$T\partial M$. We assume for simplicity that $g^{TM}$ is of {\it product structure} near
the boundary $\partial M$.
Let $S(TX)=S_+(TX)\oplus S_-(TX)$ be the ${\bf Z}_2$-graded Hermitian vector
bundle of $(TX, g^{TX})$-spinors.

Since now $M$ has a nonempty boundary $\partial M$, the associated Dirac operator
$D_+:\Gamma(S_+(TM))\rightarrow \Gamma(S_-(TM))$ is {\it not} elliptic. To get
an elliptic problem, one needs to introduce an elliptic boundary condition
for $D_+$, and this was achieved by Atiyah, Patodi and Singer in [APS]. It is remarkable
that this boundary condition, to be described right now, is {\it global} in nature.

First of all, the Dirac operator $D_+$ induces canonically a formally self-adjoint
first order elliptic differential operator
$$D_{\partial M} :\Gamma\left(S_+(TM)|_{\partial M}\right)\rightarrow
 \Gamma\left(S_+(TM)|_{\partial M}\right),$$
which is called the induced Dirac operator on the boundary $\partial M$.

Clearly, the $L^2$-completion of
$S_+ (TM)|_{\partial M} $ admits an orthogonal decomposition
$$L^2\left(S_+ (TM)|_{\partial X}\right) =
\bigoplus_{\lambda \in  {\rm Spec}(D_{\partial M} )} E_\lambda  ,\eqno (3.1)$$
where $E_\lambda$ is the eigenspace of $\lambda$.

Let  $L^2_{\geq 0} (S_+ (TM)|_{\partial M}) $
denote the direct sum of the eigenspaces $E_\lambda$ associated to the
eigenvalues $\lambda \geq 0$.
Let $P_{\geq 0} $ denote the orthogonal projection from
$L^2(S_+ (TM)|_{\partial M}) $ to
$L^2_{\geq 0}(S_+ (TM)|_{\partial M}) $.
We  call $P_{\geq 0} $ the {\it Atiyah-Patodi-Singer projection}
associated to $D_{\partial M} $, to emphasize its role in [APS].

Then by [APS], the  boundary  problem
$$\left(D_+, P_{\geq 0} \right) :\left\{ u:u\in \Gamma (S_+(TM)) ,\
P_{\geq 0}\left( u|_{\partial M}\right)=0\right\}\rightarrow \Gamma (S_-(TM)),
\eqno (3.2)$$
is Fredholm. We call this elliptic boundary  problem
the Atiyah-Patodi-Singer boundary problem associated
to $D_+$. We denote by $ {\rm ind}\, (D_+, P_{\geq 0}) $ the index of
the Fredholm operator (3.2).

$\ $

\noindent {\bf The Atiyah-Patodi-Singer index theorem} {\it The following identity holds,}
$${\rm ind}\, \left(D_+, P_{\geq 0}\right) =\int_M
\widehat{A} \left( {R^{TM}\over 2\pi}\right) -\overline{\eta}\left(D_{\partial M}
\right).\eqno(3.3)$$

$\ $

The boundary correction term $\overline{\eta}(D_{\partial M})$ appearing in the
 right hand side of (3.3)
is a spectral invariant associated to the induced Dirac operator $D_{\partial M}$
on $\partial M$. It is defined as follows: for any complex number $s\in {\bf C}$
with ${\rm Re}(s)>\dim M$, define
$$\eta\left(D_{\partial M},s\right) =\sum_{\lambda \in  {\rm Spec}(D_{\partial M} )}
{{\rm sgn}(\lambda)\over |\lambda|^s}.\eqno(3.4)$$
By using the heat kernel method, one can show easily that $\eta(D_{\partial M} ,s) $ can be
extended to a meromorphic function on ${\bf C}$, which is holomorphic at $s=0$.
Following [APS], we then define
$$\overline{\eta}\left(D_{\partial M}\right)=
{\dim\left( \ker D_{\partial M}\right)
 + \eta\left(D_{\partial M} ,0\right) \over 2}\eqno(3.5)$$
and call it the (reduced) eta invariant of $D_{\partial M} $.

The eta invariants of Dirac operators have played
important roles in many aspects of topology, geometry and mathematical physics.

In the next sections, we will discuss the role of
eta invariants in the heat kernel approaches to the index theorems on odd
dimensional manifolds.

\section{Heat kernels and the index theorem on odd dimensional manifolds}

\vskip-5mm \hspace{5mm}

Let $M$ be now an {\it odd} dimensional smooth closed oriented {\it spin} manifold.
Let $g^{TM}$ be a Riemannian metric on $TM$ and $S(TM)$ the associated Hermitian
vector bundle of $(TM,g^{TM})$-spinors.\footnote{Since now $M$ is of odd dimension,
the bundle of spinors does not admit a ${\bf Z}_2$-graded structure.} In this case,
the associated {\it Dirac} operator $D:\Gamma(TM)\rightarrow \Gamma(TM)$ is (formally)
{\it self-adjoint}.\footnote{In fact, if $M$ bounds an even dimensional spin manifold,
then $D$ can be thought of as the induced Dirac operator on boundary appearing in the
previous section.} Thus, one can proceed as in Section 3 to construct the
Atiyah-Patodi-Singer projection
$$P_{\geq 0} :L^2(S (TM))\rightarrow
L^2_{\geq 0}(S (TM)) .$$

Now consider the trivial vector bundle ${\bf C}^N$ over $M$.
We equip ${\bf C}^N$ with the canonical trivial metric and connection.
Then $P_{\geq 0}$ extends naturally
to an orthogonal projection from $L^2(S(TM)\otimes {\bf C}^N)$ to
$L^2_{\geq 0}(S(TM)\otimes {\bf C}^N)$ by acting
as identity on ${\bf C}^N$. We still denote
this extension by $P_{\geq 0}$.

On the other hand, let
$$g:M\rightarrow U(N)$$
be a smooth map from $M$ to the unitary group $U(N)$. Then
$g$ can be interpreted as automorphism of the trivial complex vector bundle ${\bf C}^N$.
 Moreover $g$
extends naturally to an action on $L^2(S(TM)\otimes {\bf C}^N)$ by acting as identity on
$L^2(S(TM))$. We still denote this extended action by $g$.

With the above data given, one can define a  {\it Toeplitz} operator $T_g$ as follows,
$$T_g=P_{\geq 0}gP_{\geq 0}:L^2_{\geq 0}\left(S(TM)\otimes {\bf C}^N\right)
\longrightarrow L^2_{\geq 0}\left( S(TM)\otimes {\bf C}^N\right).\eqno (4.1)$$

The first important fact is that $T_g$ is a Fredholm operator. Moreover,
it is equivalent to an elliptic
pseudodifferential operator of order zero. Thus one can compute its index by using the
Atiyah-Singer index theorem [AS2], as was indicated in the paper of Baum and Douglas [BD], and
the result is
$${\rm ind}\, T_g=-\left\langle \widehat{A}(TM) {\rm ch}(g),[M]\right\rangle ,
\eqno(4.2)$$
where ${\rm ch}(g)$ is the odd
Chern character associated to $g$.

There is also an analytic proof of (4.2) by using  heat kernels. For this one first
applies a result of
Booss and Wojciechowski (cf. [BW]) to show that the computation of
${\rm ind}\, T_g$ is equivalent
to the computation of the spectral flow of the
linear family of self-adjoint elliptic operators,
acting of $\Gamma(S(TM)\otimes {\bf C}^N)$, which connects $D$ and $gDg^{-1} $.
The resulting spectral flow can then be computed by  variations of $\eta$-invariants,
where the
heat kernels are naturally involved.

The above ideas have been extended in [DZ1] to give a heat kernel proof of a family
extension of (4.2).

\section{An index theorem for Toeplitz operators on odd dimensional manifolds with
boundary}

\vskip-5mm \hspace{5mm}

In this section, we describe an extension of (4.2) to the case of manifolds with boundary,
which was proved recently in my paper with Xianzhe Dai [DZ2]. This result can be thought of
as an odd dimensional analogue of the Atiyah-Patodi-Singer index theorem described
in Section 3.

This section is divided into three subsections. In Subsection 4.1,
we extend the definition of
 Toeplitz operators to the case of manifolds with boundary. In  Subsection 4.2,
we define
an $\eta$-invariant for cylinders which will appear in the statement of the main result
to be described in Subsection 4.3.

\subsection{Toeplitz operators on manifolds with boundary}\vskip-5mm \hspace{5mm}

Let $M$ be an odd dimensional oriented {\it spin} manifold with (nonempty)
boundary $\partial M$.
Then $\partial M$ is also oriented and spin.
Let $g^{TM}$ be a Riemannian metric on $TM$ such that it is of
product structure near the boundary
$\partial M$. Let $S(TM)$ be the Hermitian
bundle of spinors associated to $(M,g^{TM})$.
Since $\partial M\neq \emptyset$,
 the {\it Dirac} operator  $D:\Gamma(S(TM))\rightarrow \Gamma(S(TM))$ is no longer
elliptic. To get an elliptic operator, one needs to impose suitable boundary
conditions, and it turns out that again we will adopt the boundary
conditions introduced by Atiyah, Patodi and Singer [APS].

Let $D_{\partial M} :\Gamma(S(TM)|_{\partial M} )\rightarrow \Gamma(S(TM)|_{\partial M} )$
be the canonically induced Dirac operator on the boundary $\partial M$.
Then $D_{\partial M} $ is elliptic and (formally) self-adjoint. For simplicity,
we assume here that $D_{\partial M} $ is {\it invertible}, that is, $\ker D_{\partial M} =0$.

Let $P_{\partial M,\geq 0}$ denote the Atiyah-Patodi-Singer  projection from
$L^2(S(TM)|_{\partial M})$ to $L^2_{\geq 0}(S(TM)|_{\partial M})$.
Then $(D, P_{\partial M,\geq 0} )$ forms
a {\it self-adjoint} elliptic boundary
problem.
We will also denote the corresponding elliptic self-adjoint operator by
$D_{P_{\partial M,\geq 0}}$.

Let $L^2_{P_{\partial M,\geq 0},\geq 0}(S(TM))$ be the space of the direct sum
of eigenspaces of non-negative
eigenvalues of $D_{P_{\partial M,\geq 0}}$.
Let $P_{P_{\partial M,\geq 0},\geq 0}$ denote the orthogonal projection from
$L^2(S(TM))$ to $L^2_{P_{\partial M,\geq 0},\geq 0}(S(TM))$.

Now let ${\bf C}^N$ be the trivial complex vector bundle over
$M$ of rank $N$, which carries the trivial Hermitian metric and the trivial Hermitian connection.
We extend $P_{P_{\partial M,\geq 0},\geq 0}$ to act as identity on ${\bf C}^N$.

Let $g:M\rightarrow U(N)$ be a smooth unitary automorphism of ${\bf C}^N$.
Then $g$ extends to an action on $S(TM)\otimes {\bf C}^N$ by acting as identity
on $S(TM)$.

Since $g$ is unitary, one verifies easily that the operator
$gP_{\partial M,\geq 0} g^{-1}$
is an orthogonal projection on
$L^2((S(TM)\otimes {\bf C}^N)|_{\partial M})$, and that
$gP_{\partial M,\geq 0} g^{-1}-
P_{\partial M,\geq 0} $ is a pseudodifferential operator
of order less than zero. Moreover, the pair
$(D,gP_{\partial M,\geq 0} g^{-1})$ forms
a {\it self-adjoint} elliptic boundary problem.
We denote its associated elliptic self-adjoint operator
by $D_{gP_{\partial M,\geq 0} g^{-1}}$.

Let $L^2_{gP_{\partial M,\geq 0} g^{-1},\geq 0}
(S(TM)\otimes {\bf C}^N)$ be the space of the direct sum
of eigenspaces of nonnegative
eigenvalues of $D_{gP_{\partial M,\geq 0}g^{-1}}$.
Let $P_{gP_{\partial M,\geq 0}g^{-1},\geq 0}$
denote the orthogonal projection   from $L^2(S(TM)\otimes {\bf C}^N)$ to
$L^2_{gP_{\partial M,\geq 0}g^{-1},\geq 0}
(S(TM)\otimes {\bf C}^N)$.

Clearly, if $s\in L^2(S(TM)\otimes {\bf C}^N)$ verifies
$P_{\partial M,\geq 0}(s|_{\partial M})=0$, then
$gs$ verifies
$$gP_{\partial M,\geq 0}g^{-1}\left((gs)|_{\partial M}\right)=0.$$

$\ $

\noindent {\bf Definition 5.1} {The {\it Toeplitz operator} $T_g$ is defined by
$$T_g=P_{gP_{\partial M,\geq 0}g^{-1},\geq 0}
g P_{P_{\partial M,\geq 0} ,\geq 0}:$$
$$L^2_{P_{\partial M,\geq 0},\geq 0}
\left(S(TM)\otimes {\bf C}^N\right) \rightarrow
L^2_{gP_{\partial M,\geq 0}g^{-1},\geq 0}
\left(S(TM)\otimes {\bf C}^N\right). $$}

One verifies that $T_g$ is a Fredholm operator. The main result of this section
 evaluates  the  index of $T_g$  by more geometric quantities.

\subsection{An {\boldmath $\eta$}-invariant  associated to {\boldmath $g$}}
\vskip-5mm \hspace{5mm}

We consider the cylinder $[0,1]\times \partial M$. Clearly, the restriction of $g$
on $\partial M$ extends canonically to this cylinder.

Let $D|_{[0,1]\times \partial M}$ be the restriction of
$D$ on $[0,1]\times \partial M$. We equip the boundary condition
$P_{\partial M,\geq 0}$ at $\{0\}\times \partial M$
and the boundary condition ${\rm Id}-gP_{\partial M,\geq 0}g^{-1}$
at $\{1\}\times \partial M$. Then
$(D|_{[0,1]\times \partial M}, P_{\partial M,\geq 0},
{\rm Id}-gP_{\partial M,\geq 0}g^{-1})$ forms a self-adjoint
elliptic boundary problem. We denote the corresponding elliptic self-adjoint operator
by $D_{P_{\partial M,\geq 0},
gP_{\partial M,\geq 0}g^{-1}}$.

Let $\eta(D_{P_{\partial M,\geq 0},
gP_{\partial M,\geq 0}g^{-1}},s)$ be the $\eta$-function of
$s\in {\bf C}$ which, when ${\rm Re}(s)>>0$, is defined by
$$\eta\left(D_{P_{\partial M,\geq 0},
gP_{\partial M,\geq 0}g^{-1}},s\right)
=\sum_{\lambda \neq 0}{{\rm sgn}(\lambda)\over |\lambda|^s},$$
where $\lambda$ runs through the nonzero eigenvalues of
$D_{P_{\partial M,\geq 0},
gP_{\partial M,\geq 0}g^{-1}}$.

It is proved in [DZ2] that under our situation, $\eta(D_{P_{\partial M,\geq 0},
gP_{\partial M,\geq 0}g^{-1}},s)$ can be extended to a meromorphic function on
${\bf C}$ which is holomorphic at $s=0$.

Let $\overline{\eta}(D_{P_{\partial M,\geq 0},
gP_{\partial M,\geq 0}g^{-1}})$ be the reduced $\eta$-invariant defined by
$$\overline{\eta}\left(D_{P_{\partial M,\geq 0},
gP_{\partial M,\geq 0}g^{-1}}\right)$$
$$={\dim \ker \left(D_{P_{\partial M,\geq 0},
gP_{\partial M,\geq 0}g^{-1}}\right)+\eta\left(D_{P_{\partial M,\geq 0},
gP_{\partial M,\geq 0}g^{-1}}\right)\over 2}.$$

\subsection{An index theorem for {\boldmath $T_g$}}\vskip-5mm \hspace{5mm}

Let $\nabla^{TM}$ be the Levi-Civita connection associated to the Riemannian metric
$g^{TM}$. Let $R^{TM}=(\nabla^{TM})^2$
be the curvature of $\nabla^{TM}$.
Also, we use $d$ to denote the trivial connection on the trivial vector bundle ${\bf C}^N$
over $M$. Then $g^{-1}dg$ is a $\Gamma({\rm End}({\bf C}^N))$ valued 1-form
over $M$.

Let ${\rm ch}(g,d)$ denote the odd Chern character form (cf. [Z]) of $(g, d)$ defined by
$${\rm ch}(g,d)=
\sum_{n=0}^{(\dim M-1)/2} {n!\over (2n+1)!}
\left({1\over 2\pi\sqrt{-1}}\right)^{n+1}{\rm Tr}\left[\left(g^{-1}dg\right)^{2n+1}\right].
$$

Let ${\cal P}_M$ denote the Calder\'on projection associated to $D$ on $M$
(cf. [BW]). Then ${\cal P}_M$ is an orthogonal projection on
$L^2((S(TM)\otimes {\bf C}^N)|_{\partial M})$,
and that ${\cal P}_M -P_{\partial M,\geq 0}$ is a pseudodifferential operator of
order less than zero.

Let $\tau_\mu (P_{\partial M,\geq 0}, gP_{\partial M,\geq 0}g^{-1}, {\cal P}_M)\in {\bf Z}$
be the Maslov triple index in the sense of
Kirk and Lesch [KL, Definition 6.8].

We can now state the main result of [DZ2], which generalizes an old result of
Douglas and Wojciechowski [DoW], as follows.

$\ $

\noindent {\bf Theorem 5.2} {\it The following identity holds,
$${\rm ind}\, T_g=-\int_M
\widehat{A}\left({R^{TM}\over 2\pi }\right)
{\rm ch}(g,d)
+\overline{\eta}\left(D_{P_{\partial M,\geq 0},
gP_{\partial M,\geq 0}g^{-1}}\right) $$
$$-\tau_\mu \left(P_{\partial M,\geq 0}, gP_{\partial M,\geq 0}g^{-1},
{\cal P}_M\right) .$$}

The following immediate consequence is of independent interests.

$\ $

\noindent {\bf Corollary 5.3} {\it The number
$$\int_M
\widehat{A}\left({R^{TM}\over 2\pi }\right)
{\rm ch}(g,d)-\overline{\eta}\left(D_{P_{\partial M,\geq 0},
gP_{\partial M,\geq 0}g^{-1}}\right) $$
is an integer.}

$\ $

The strategy of the proof of Theorem 5.2 given in [DZ2] is the same as
that of the heat kernel proof of
(4.2). However, due to the appearance of the boundary $\partial M$,
one encounters new difficulties. To overcome these difficulties,
one makes use of the recent result on the splittings of $\eta$ invariants
(cf. [KL]) as well as some ideas involved in the
Connes-Moscovici local index theorem in noncommutative geometry [CM]
(see also [CH]). Moreover, the local index calculations appearing near $\partial M$
is highly nontrivial. We refer to [DZ2] for more details.

\label{lastpage}

\end{document}